\documentclass[a4paper]{amsart}

\title{BPS  state counting on singular varieties}
\author{Elizabeth \ Gasparim}
\address{E. Gasparim\\ IMECC-UNICAMP, Cidade Universitária "Zeferino Vaz",
Campinas, SP, 13083-859, Brasil.}
\email{etgasparim@gmail.com}
\author{ Thomas \ K\"{o}ppe}
\address{T. K\"oppe, Department of Mathematics, King's College London,
Strand  WC2R 2LS, UK.}
\email{thomas.koeppe@kcl.ac.uk}
\author{ Pushan \ Majumdar}
\address{P. Majumdar\\Dept. of Theoretical Physics,
Indian Association for the Cultivation of Science, Calcutta 700 032, India.}
\email{tppm@iacs.res.in}
\author{ Koushik \ Ray}
\address{ K. Ray\\Dept. of Theoretical Physics,
Indian Association for the Cultivation of Science, Calcutta 700 032, India.}
\email{koushik@iacs.res.in}
\thanks{Generous support of the Royal Society
under Grant 2008/R2 for joint research projects between the UK and India
made this collaboration possible.} 
\date{}
\usepackage[margin=3cm]{geometry}
\usepackage[dvipsnames]{xcolor}
\usepackage[pdftitle={\@title}, pdfauthor={\@author}, pdfpagelayout=OneColumn, colorlinks=true, linkcolor=black, urlcolor=black, citecolor=black]{hyperref}
\usepackage{amsmath,amsthm,amssymb,amscd,bm,cancel,booktabs,tikz}
\usepackage{pstricks,pst-3dplot,pst-slpe,pst-grad,
pst-fill,pst-plot,pst-tree,pst-3d}
\usepackage{enumerate}

\numberwithin{equation}{section}

\DeclareMathOperator{\ev}{ev}                        
\DeclareMathOperator{\spec}{Spec}
\newcommand{\vd}{\text{\rm vd}}
\newcommand{\vir}{\text{\rm vir}}
\newcommand{\om}{\kern.5em\overline{\phantom{M}}\kern-1.2em\mathcal{M}}

\newcommand{\R}{\mathbb{R}}
\newcommand{\C}{\mathbb{C}}
\newcommand{\Z}{\mathbb{Z}}
\newcommand{\ce}{\mathrel{\mathop:}=}                
\newcommand{\ec}{=\mathrel{\mathop:}}                
\newcommand{\abs}[1]{\left\lvert#1\right\rvert}      
\newcommand{\oton}[2][n]{#2_1, \dotsc, #2_{#1}}      
\newcommand{\eqand}[1][and]{\qquad\text{#1}\qquad}   
\newcommand{\cp}{\mathbb{P}}
\newcommand{\xcy}{X_{\text{CY}}}
\newcommand{\ns}{\pscircle[linestyle=none,fillstyle=ccslope,
slopeend=black!60!white,slopebegin=black,slopecenter=.54 .46](0,0){.14 }}
\newcommand{\nsp}{\pscircle[linestyle=none,fillstyle=ccslope,
slopeend=black,slopebegin=white,slopecenter=.54 .46](0,0){.14}}

\newtheorem{theorem}{Theorem}[section]
\newtheorem*{theorem*}{Theorem}

\newtheorem*{lemma*}{Lemma}
\newtheorem{proposition}[theorem]{Proposition}
\newtheorem*{proposition*}{Proposition}
\newtheorem{corollary}[theorem]{Corollary}
\newtheorem*{corollary*}{Corollary}

\newtheorem*{conjecture*}{Conjecture}
\theoremstyle{definition}
\newtheorem{definition}[theorem]{Definition}
\newtheorem*{definition*}{Definition}
\newtheorem{remark}[theorem]{Remark}
\newtheorem*{remark*}{Remark}
\newtheorem{example}[theorem]{Example}
\newtheorem*{example*}{Example}

\begin{document}

\begin{abstract}
We define new partition functions for theories with targets 
on toric singularities via
products of old partition functions on  crepant resolutions.
We compute explicit examples and
show that the  new partition functions
turn out to be homogeneous on MacMahon factors. 
\end{abstract}
\maketitle
\tableofcontents
\section{Motivation for counting BPS states}
BPS states are minimal energy states of supersymmetric field theories. 
These special states have had a crucial role in establishing various
duality symmetries of Superstring theory. One of the reasons for their
pivotal role in studying dualities stems from the availability of
information on exact masses and degeneracies of these states.
The counting of BPS states is of great interest to string theory and
supergravity. In certain instances 
the counting of BPS states agrees with the
counting of extremal black holes
\cite{S1, IS}. In some cases 
the string partition function matches with the black
hole partition function, leading to a precise equivalence between the 
black hole entropy and the
statistical entropy associated with an ensemble of BPS states \cite{S2}.
Degeneracy of states is encoded in a partition function.  
Degeneracy of BPS D-branes in string theory 
depends on the background geometry. 
The spectrum of BPS D-branes changes across walls in the moduli space. 
As the moduli of the background is
varied, the number of states can jump across  walls of marginal
stability. The walls thus partition the moduli space into chambers.
In other words, across wall a BPS state may
disappear, or `decay', giving rise to a different spectrum of BPS
states. The original BPS state is thus stable in a specific chamber, 
while the decay products are stable in another. Indeed, when
D-branes are realized as BPS states, they are defined by the stable
BPS states only.
Characterising the jumps of degeneracy of BPS states
across walls in the moduli space, notwithstanding the continuity of
appropriate correlation functions, has been of immense interest
recently \cite{OSY,Na,KS1,KS2,MMNS}.
These studies unearthed a rich mathematical structure within
the scope of topological string theories.

A class of BPS states in topological string theories is furnished by
D-branes wrapping homology cycles of the target space. These D-branes
as well as their bound states are described as objects in the derived
category of coherent sheaves of the target space or objects in the
Fukaya category, within the scope of the topological B or A models,
respectively. On a Calabi-Yau target the walls of marginal stability
are detected from the alignment of charges of the D-branes in the
spectrum. Across a wall a D-brane decays into a finite or infinite
collection of branes, with the charge of the parent brane aligning
with the totality of charges of the products on the wall. The partition 
function of these branes can be calculated giving their degeneracies.

The partition function of the A-model generates the
Gromov--Witten (GW) invariants of Calabi-Yau threefolds from the
world-sheet perspective. From the target space perspective, it counts
the Gopakumar--Vafa (GV) invariants.  The GW invariants count holomorphic
curves on the threefold, whereas the GV invariants count BPS states of
spinning black holes in 5 dimensions obtained from M2-branes in
M-theory on the Calabi-Yau threefold \cite{AOVY}. Considering the
topological A-model on the target
$\R^3\times X\times S^1$, 
where $X$ denotes the Calabi-Yau space without four-cycles and $S^1$
designates the compact Euclidean temporal direction, the partition
function also counts the number of D0- and D2-brane bound states on a
single D6-brane wrapped on $X$. M5-branes wrapping
four-cycles in $X$ may form bound states with M2-branes; these
complications do not arise in the absence of four-cycles in $X$
\cite{AOVY}.  From another point of view
the partition function of the A-model is
also the generating function of
the Donaldson--Thomas invariants in appropriate variables. 
Thus the study of the degeneracy of
states relates the GW, GV and DT invariants.

For a singular variety, for example an orbifold, 
the product of the partition functions for all its 
crepant resolutions may be considered. 
The homology groups of the crepant resolutions are isomorphic.
For the BPS D-branes the crepant resolutions 
correspond to different spectra of stable objects in different chambers 
with the isomorphism of homologies given by Seiberg duality.
The product partition function then corresponds to a quiver variety, 
which is realised near the singularity or the orbifold point, possessing 
a derived equivalence with the crepant 
resolutions \cite{Sz,Na,Y}. However, different isomorphisms of homologies yield 
different partition functions. Here we define a partition function for 
the generalised conifolds as the product of
the crepant resolutions as above, but the isomorphism of the second 
homology groups is given by a direct identification of elements in terms of 
certain formal variables under a canonical ordering.
In proving the main theorem on the homogeneity of the new partition function
we use a probabilistic argument which appears to
relate the exponent of homogeneity to some kind of degeneracy of the singular
variety. Finally, we discuss some combinatoric aspects of the T-dual
type-IIA brane configurations with NS and NS' branes
corresponding to the crepant resolutions of
$C_{m,n}$, which is related to the partition function of the quiver
variety. We write down explicit formulas for the generalised conifold
$C_{1,3}$ and compare the two partition functions.
\section{New partition function via formal identification 
and main results}
Let $X$ be a singular variety admitting a finite
collection of crepant resolutions $X^t \to X$ for an index $t \in
\mathcal{T}$, $\abs{\mathcal{T}} < \infty$. 
If a singular variety admits crepant resolutions each of which have trivial
canonical bundle, then it will
be called a singular Calabi-Yau variety. Let $X$ be a singular Calabi-Yau
variety. Let us further assume that
a partition function $Z_\text{old}(Y; Q, \dotsc)$ is defined
for a smooth Calabi-Yau space $Y$, where $Q = (Q_1, Q_2, \dotsc)$
are formal variables corresponding to a basis of $H_2(Y;\mathbb{Z})$.
Finally, we suppose that $H_2(X^{s};\mathbb{Z}) \cong H_2(X^{t};\mathbb{Z})$
for all $s,t \in \mathcal{T}$.

We then define a new partition function as the product of partition functions
of the resolutions,
\[ Z_\text{new}(X;Q,\dotsc) \ce \prod_{t \in \mathcal{T}} Z_\text{old}(X^t; Q^t, \dotsc). \]
The new partition function contains information about all crepant
resolutions of $X$ and may thus be regarded as pertaining to the singular
space $X$ itself. In the product we do not include partial resolutions 
as they are contained in the full resolutions and their inclusion will but
cause non-illuminating repetitions. 
This approach can be applied to various  partition
functions defined for Calabi-Yau spaces. 
In this paper, we restrict to partition functions of curve-counting type 
such as the Gromov--Witten and
the Donaldson--Thomas partition functions. 

The properties of the new partition function depend on the prescribed 
isomorphism of second homologies of the resolutions. 
Assuming a canonical ordering of elements of $H_2(X^t,\Z)$, for all $t\in\mathcal{T}$,
we identify the formal variables $Q$ among all the
resolutions giving the isomorphism of homologies by setting
\begin{equation}\label{id} Q_i^s = Q_i^t \ec Q_i \text{ \ for all \ } s,t 
\in \mathcal{T}.\end{equation}
Presence of four-cycles in the resolutions complicates the ordering of second 
homologies. We shall restrict to varieties whose crepant resolutions do not possess 
homology four-cycles. 

By a \emph{singular Calabi-Yau threefold without contractible curves 
and/or compact $4$-cycles} we refer to a singular Calabi-Yau variety
admitting  crepant resolutions, the latter containing no
contractible curve and/or compact $4$-cycle. 
We  prove the following 
\begin{theorem*}
Let $X$ be a singular toric Calabi-Yau threefold
defined as a subset of $\C^4$ by 
$X=\C[x,y,z,w]/\langle xy-z^mw^n\rangle$, where $m$ and $n$ are
integers. 
Let $Z(Y;q,Q)$ be a partition
function of curve-counting type \text{\rm (Definition~\ref{cc})}.
Then the total partition function
\[ Z_\text{\rm tot}(X;q,Q) \ce 
\prod_{\substack{Y\\Y \to X}} Z(Y;q,Q),\]
where the product ranges over all crepant resolutions of $X$, is homogeneous
\text{\rm (Definition~\ref{def.hom})}
of degree 
\[  d = \frac{(m^2-m+n^2-n-2mn)(m+n-2)!}{m!n!}. \]
\end{theorem*}
In performing  curve counting the Calabi-Yau
space is allowed to have contractible curves as well
(Corollary \ref{contrac}) in particular obtaining a counting of 
BPS states via the topological string partition 
function (Corollary \ref{BH}).
\section{The mathematics of curve counting}
\subsection{Gromov--Witten theory}
\begin{definition}
By a \emph{curve} we mean a reduced scheme $C$ of pure dimension one.
The \emph{genus} of $C$ is $g(C) \ce h^1(C;\;\mathcal{O}_C)$.
\end{definition}

\begin{corollary}
A connected curve $C$ of genus $0$ is a tree of rational curves.
\end{corollary}

\begin{definition}
An $n$-pointed curve $\bigl(C;\oton{P}\bigr)$ is called
\emph{prestable} if every point of $C$ is either smooth or a node
singularity and the points $\oton{P}$ are smooth. A map $f \colon C
\to X$ is called \emph{stable} if $\bigl(C;\oton{P}\bigr)$ is
prestable and there are at least three marked or singular points on
each contracted component.
\end{definition}

\begin{remark}
Stability prohibits first-order infinitesimal deformations to the map $f$.
\end{remark}

Let us denote by $\om_{g,n}(X, \beta)$ the collection of maps from stable,
$n$-pointed curves of genus $g$ into $X$ for which
\[ [f(C)] = f_*[C] = \beta \in H_2(X;\;\mathbb{Z}) \text{ .} \]
Behrend and Fantechi \cite{BF1} showed that this has a coarse moduli
(Deligne--Mumford) stack, Vistoli \cite{V} studied the intersection theory on
$\om_{g,n}(X, \beta)$ and constructed a perfect obstruction theory,
and \cite{BF1} showed that there exists a virtual fundamental
class of virtual dimension
\[ \vd = (1-g) (\dim X - 3) - K_X(\beta) + n \text{ .} \]
(We assume that $X$ does in fact have a canonical class $K_X
\in H^2(X;\;\mathbb{Z})$, e.g.\ if $X$ is smooth.) Consequently,
dimension of the classes $[\om_{g,n}(X,\beta)]^\vir$ is independent of
$\beta$ when $K_X=0$, that is, when $X$ is Calabi-Yau. Moreover, the
unpointed moduli $\om_{0,0}(X,\beta)$ has virtual dimension zero for
all $g$ if $\dim X = 3$, so on a three-dimensional Calabi-Yau,
$\om_{0,0}(X,\beta)$ really ``counts curves''.
\begin{definition}
Assume that $g(C) = 0$. Let
\[ \ev_i \colon \om_{0,n}(X,\beta) \to X \text{ ,}\quad
   \bigl( f \colon (C; \oton{P}) \to X \bigr) \mapsto f(P_i) \]
be the $i^\text{th}$ evaluation map. Assume that $\sum_{i=1}^n
\deg(\gamma_i) = \vd$ for some $\gamma_i \in H^*(\om_{0,n}(X,
\beta))$. Then the \emph{genus-$0$ Gromov--Witten invariants} are
\[ \langle \oton\gamma\rangle_\beta \ce \ev_1^*(\gamma_1) \cup \dotsb \cup
   \ev_n^*(\gamma_n) \cap [\om_{0,n}(X,\beta)]^\vir \text{ .} \]
For higher genera, the definition of the Gromov--Witten invariants
requires the introduction of additional data, called \emph{descendent
fields}. Since we require only genus $0$ for our purposes, we
refer the interested reader to \cite[\S\,2]{MNOP2}.
\end{definition}

When $\dim X = 3$, $X$ is Calabi-Yau (i.e.\ $K_X=0)$, arbitrary genus
$g$ and $n=0$, we have the \emph{unmarked Gromov--Witten invariants}
\[ N_{g,\beta}(X) \ce \int_{[\kern.3em\overline{\phantom{N}}\kern-.85em\mathcal{M}_{g,0}(X,\beta)]^\vir} 1. \]
\begin{example}
If $X=\{\text{pt}\}$, then $\om_{g,n}(X,\beta) = \om_{g,n}$, the
moduli of $n$-pointed curves.
\end{example}
\begin{example}
For $X = \cp^1$, the genus-$0$ Gromov--Witten 
invariants are just the Hurwitz numbers.
\end{example}
The (unmarked) Gromov--Witten invariants are usually assembled into an
unreduced and a reduced generating function, respectively
\begin{equation*} 
F(X;u,v) = \sum_{\beta} \sum_{g\geq0} N_{g,\beta}(X) u^{2g-2} v^\beta,
\end{equation*} 
and
\begin{equation*} 
F'(X;u,v) = \sum_{\beta\neq0} \sum_{g\geq0} N_{g,\beta}(X) u^{2g-2} v^\beta,
\end{equation*} 
where $v=(\oton[r]v)$ is an appropriate vector that can be paired with
the $r$ generators of $H_2(X;\;\mathbb{Z})$. The unreduced and reduced
\emph{Gromov--Witten partition functions} are, respectively,
\begin{equation*} 
Z_{\text{GW}}(X;u,v)= 
\exp F(X;u,v) = 1 + \sum_{\beta} Z_{\text{GW}}(X;u)_\beta v^\beta
\end{equation*} 
and
\begin{equation*} 
Z'_{\text{GW}}(X;u,v) = \exp F'(X;u,v) = 1 + \sum_{\beta \neq 0}
Z'_{\text{GW}}(X;u)_\beta v^\beta,
\end{equation*} 
where the last expressions define the homogeneous terms $Z(X;u)_\beta$
and $Z'(X;u)_\beta$ of degree $\beta$.
\subsection{Donaldson--Thomas theory}
An \emph{ideal subsheaf} of $\mathcal{O}_X$ is a sheaf $\mathcal{I}$
such that $\mathcal{I}(U)$ is an ideal in $\mathcal{O}_X(U)$ for each
open set $U \subseteq X$. Alternatively, it is a torsion-free rank-$1$
sheaf with trivial determinant. It follows that $\mathcal{I}^{\vee\vee}
\cong \mathcal{O}_X$. Thus the evaluation map determines a quotient
\begin{equation}\label{eq.ideal}
  0 \longrightarrow \mathcal{I} \xrightarrow{\ \ev\ } \mathcal{I}^{\vee\vee}
  \cong \mathcal{O}_X \longrightarrow \mathcal{O}_X\bigl/\mathcal{IO}_X = \imath_*\mathcal{O}_Y
  \longrightarrow 0 \text{ ,}
\end{equation}
where $Y \subseteq X$ is the support of the quotient and
$\mathcal{O}_Y \ce (\mathcal{O}_X\bigl/\mathcal{IO}_X)\rvert_Y$ is the
structure sheaf of the corresponding subspace. Let $[Y] \in H_2(X;\;\mathbb{Z})$
denote the cycle class determined by the $1$-dimensional components of $Y$ with
their intrinsic multiplicities. We denote by
\[ I_n(X,\beta) \]
the Hilbert scheme of ideal sheaves $\mathcal{I} \subset \mathcal{O}_X$
for which the quotient $Y$ in \eqref{eq.ideal} has dimension at most
$1$, $\chi(\mathcal{O}_Y) = n$ and $[Y] = \beta \in H_2(X; \;
\mathbb{Z})$.

The work of Donaldson and Thomas was to show that $I_n(X,\beta)$ has a
canonical perfect obstruction theory (originally when $X$ is smooth,
projective and $-K_X$ has non-zero sections) and a virtual fundamental
class $[I_n(X,\beta)]^\vir$ of virtual dimension $\int_\beta c_1(T_X)
= -K_X(\beta)$. If $X$ is a smooth, projective Calabi-Yau threefold,
then the virtual dimension is zero, and we write
\[ \widetilde{N}_{n,\beta}(X) \ce \int_{[I_n(X,\beta)]^\vir} 1 \]
for the ``number'' of such ideal sheaves. We assemble these
numbers into an (unreduced) partition function,
\[ Z_{\text{DT}}(X;q,v) = \sum_{\beta} \sum_{n\in\mathbb{Z}} \widetilde{N}_{n,
   \beta}(X) q^n v^\beta = \sum_\beta Z_{\text{DT}}(X;q)_\beta v^\beta \text{ ,} \]
where again the last expression defines the unreduced terms of degree $\beta$. The degree-$0$ term
\[ Z_{\text{DT}}(X;q)_0 = \sum_{n\geq0} \widetilde{N}_{n, 0}(X) q^n \]
is of special importance: We define the \emph{reduced} DT partition
function as
\[ Z'_{\text{DT}}(X;q,v) = Z_{\text{DT}}(X;q,v) \bigl/ Z_{\text{DT}}(X;q)_0
   = 1 + \sum_{\beta \neq 0} Z'_{\text{DT}}(X;q)_\beta v^\beta \text{ ,} \]
once again defining the reduced terms $Z'_{\text{DT}}(X;q)_\beta$ of
degree $\beta$ implicitly.

\subsection{The MNOP Conjecture}
For a smooth  Calabi-Yau threefold $\xcy$ the MNOP conjecture 
relates the reduced GW and DT partition functions,
\[ Z'_\text{GW}(\xcy; u, v) = Z'_\text{DT}(\xcy; -e^{iu}, v), \]
signifying an equivalence between the
Gromov--Witten and Donaldson--Thomas
theories for  Calabi-Yau threefolds.
Proof of the MNOP relation was furnished 
for  toric (hence non-compact) Calabi-Yau
threefolds in \cite{MNOP1,MNOP2} and for compact Calabi-Yau
manifolds in \cite{BF2,L}.

We shall illustrate features of the new partition function using the
Donaldson--Thomas partition function, for which toric computational
techniques have been developed by \cite{LLLZ}. 
We consider a special class of threefolds admitting crepant resolutions
without compact four-cycles,
obtained as orbifolds of the conifold or their partial resolutions
\cite{mp,uranga,unge,Na}.
\section{Generalised conifolds}\label{gencon}
Given a pair of non-negative integers $m$, $n$, we 
consider the  toric varieties
\[ C_{m,n} \ce \{(x,y,z,w)| x y - z^m w^n = 0\} 
\subset \mathbb{C}^4 = \spec\mathbb{C}[x,y,z,w].\]
We suppose $n\geq m$ without any loss of generality.
Two cases arise: 
\begin{enumerate}[(i)]
\item $n>m=0$. Then $C_{0,n}$ are quotients of $\mathbb{C}^3$ by
$\Z/n\Z$ acting on a two-dimensional subspace $\mathbb{C}^2$
as  $(a,b,c) \mapsto (\varepsilon a,
\varepsilon^{-1} b, c)$, with $\varepsilon^n = 1$. These spaces
have $1$-dimensional singularities, as $C_{0,n} \cong K_n \times \mathbb{C}$,
where $K_n=\{(x,y,z)| xy-z^n =
0\}\subset\mathbb{C}^3$ is the Kleinian surface,
 with a singular point at the origin.
\item $n\geq m \geq 1$. The space
$C_{1,1}=\{(x,y,z,w)| xy-zw=0\}\subset \mathbb{C}^4$ is the \emph{conifold}. 
All other $C_{m,n}$
are obtained either as quotients of the conifold, if $n=m$, 
or through their partial resolutions, otherwise.
\end{enumerate}
\begin{figure}[h]
\begin{pspicture}(-3.4,-2)(4,4)
\psset{Alpha=35,unit=1cm}
\pstThreeDCoor[linecolor=darkgray]
\pstThreeDLine{->}(0,0,0)(1,0,0)  
\pstThreeDLine{->}(0,0,0)(0,1,0)  
\pstThreeDTriangle[fillstyle=solid,fillcolor=black!40!white,linewidth=.5pt](0,0,0)(0,0,3)(3,0,3)
\pstThreeDTriangle[fillstyle=solid,fillcolor=black!40!white,linewidth=.5pt](0,0,0)(0,0,3)(0,1,3)
\pstThreeDLine{->}(0,0,0)(-2,-2,2)    
\pstThreeDLine{->}(0,0,0)(0,-3.6,1.2) 
\pstThreeDTriangle[fillstyle=solid,fillcolor=black!60!white,linewidth=.5pt](0,0,0)(0,1,3)(2,1,3)
\pstThreeDTriangle[fillstyle=solid,fillcolor=black!80!white,linewidth=.5pt](0,0,0)(2,1,3)(3,0,3)
\pstThreeDPut(1,.5,0){$n_1$}
\pstThreeDPut(0,-3.6,1.5){$n_2$}
\pstThreeDPut(-2,-2,2.2){$n_3$}
\pstThreeDPut(-.4,1,0){$n_4$}
\pstThreeDPut(.3,0,3.25){$v_1$}
\pstThreeDPut(0,1.4,3){$v_2$}
\pstThreeDPut(1.8,1,3.25){$v_3$}
\pstThreeDPut(3.25,0,3){$v_4$}
\end{pspicture}
\begin{pspicture}(-3,-2)(4,4)
\psset{Alpha=35,unit=1cm}
\pstThreeDCoor[linecolor=darkgray]
\pstThreeDTriangle[fillstyle=solid,fillcolor=black!40!white,linewidth=.5pt](0,0,0)(0,0,3)(3,0,3)
\pstThreeDTriangle[fillstyle=solid,fillcolor=black!40,linewidth=.5pt](0,0,0)(0,0,3)(0,1,3)
%
\pstThreeDLine[linecolor=white](2,0,3)(0,1,3)
\pstThreeDLine[linecolor=white](1,0,3)(0,1,3)
\pstThreeDLine[linecolor=white](2,1,3)(2,0,3)
\pstThreeDLine[linecolor=white](1,1,3)(2,0,3)
\pstThreeDLine[linewidth=.5pt](0,0,3)(3,0,3)
%
%
\pstThreeDTriangle[fillstyle=solid,fillcolor=black!60!white,linewidth=.5pt](0,0,0)(0,1,3)(2,1,3)
\pstThreeDTriangle[fillstyle=solid,fillcolor=black!80!white,linewidth=.5pt](0,0,0)(2,1,3)(3,0,3)
\end{pspicture}
\caption{Left: The toric fan for $C_{2,3}$ and the normal vectors 
Right: A triangulation of the cone}
\label{torvar}
\end{figure}
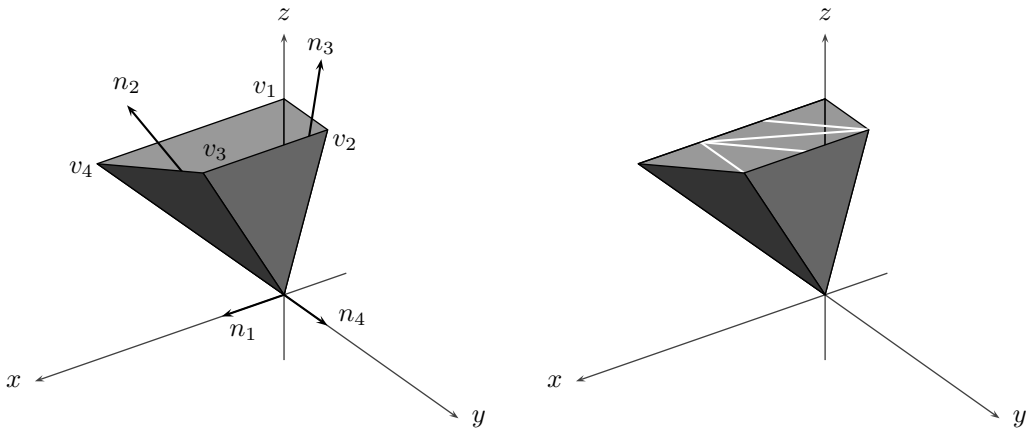
Let us briefly describe $C_{m,n}$, referred to as  a 
\emph{generalised conifold} in the sequel, as a toric variety.   
The toric fan of $C_{m,n}$ is generated by a $3$-dimensional
cone $\sigma$ with ray generators $v_i$, $i=1,2,3,4$,
which are vectors in a lattice $N$ of rank $3$ in
$\R^3$ given by the columns of the matrix
\begin{equation}
\bordermatrix{
&v_1&v_2&v_3&v_4\cr
&0&0&m&n\cr
&0&1&1&0\cr
&1&1&1&1
},
\end{equation} 
all of which lie on the height-one $z$-plane along the perimeter of 
a trapezoid, thereby rendering the canonical divisor trivial. 
The inward-pointing normals to the facets subtended by a pair of these
vectors given by their cross products, 
namely, $n_i = v_{i+1}\times v_i$, in cyclic order,
define the semigroup $S_{\sigma} = \sigma^{\vee}\cap M$, $M$ being
the dual lattice of $N$. The dual cone is 
\begin{equation}
\sigma^{\vee}= \{n\in\R^3| \langle n,v\rangle\geq 0, \forall v\in\sigma\}
\end{equation} 
The various vectors and the cone are depicted in Figure~\ref{torvar}(Left).
Then $S_{\sigma}$ is generated by the four columns
of the following matrix
\begin{equation} 
T=
\bordermatrix{
&n_1& n_2&n_3&n_4\cr
&1&0&-1&0\cr
&0&-1&m-n&1\cr
&0&1&n&0
},
\end{equation} 
which provides the toric data. The relation among these four
three-dimensional vectors is given through the kernel of $T$, 
\begin{equation}
\ker T = (1,-n,1,-m)^t. 
\end{equation} 
Hence the toric variety $C_{m,n}$ is given by the equation
\begin{equation}
x_1x_3-x_4^mx_2^n=0. 
\end{equation} 
Since all the ray generators $v_i$ lie in the height one $z$-plane,
it suffices, especially
for the purpose of exhibiting triangulations considered below, to
draw the intersection of the cone $\sigma$ with this plane.
We shall  henceforth refer to the 
trapezoidal polygon on this plane formed by
the vertices 
$\left(\begin{smallmatrix}0\\0\end{smallmatrix}\right)$, 
$\left(\begin{smallmatrix}0\\1\end{smallmatrix}\right)$, 
$\left(\begin{smallmatrix}m\\1\end{smallmatrix}\right)$, 
and $\left(\begin{smallmatrix}n\\0\end{smallmatrix}\right)$, illustrated
below, as the toric data for the variety $C_{m,n}$.
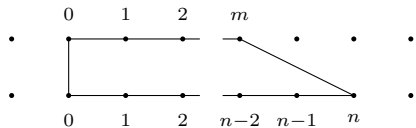
\begin{center}\begin{tikzpicture}[scale=.75]
\foreach \x in {-1,0,1,2,3,4,5,6} {
  \foreach \y in {0,1} {
    \node at (\x,\y) [circle,fill=black,scale=.2] {};
  }
}
\draw (2.3,0) 
-- (2,0) node[label=below:$\scriptstyle 2$]{} 
-- (1,0) node[label=below:$\scriptstyle 1$]{}
-- (0,0) node[label=below:$\scriptstyle 0$]{} 
-- (0,1) node[label=above:$\scriptstyle 0$]{} 
-- (1,1) node[label=above:$\scriptstyle 1$]{} 
-- (2,1) node[label=above:$\scriptstyle 2$]{} 
-- (2.3,1); 
\draw (2.7,1) -- (3,1) node[label=above:$\scriptstyle m$]{} 
-- (5,0) node[label=below:$\scriptstyle n$]{} 
-- (4,0) node[label=below:$\scriptstyle n-1$]{} 
-- (3,0) node[label=below:$\scriptstyle n-2$]{} -- (2.7,0); 
\end{tikzpicture}\end{center}
In general, blowing up the
singular locus of a generalised conifold
results in a non-Calabi-Yau variety. This can
be seen by constructing the star subdivision of the singular
subcone. The new ray generator does not lie on the
$z=1$ hyperplane. However, small resolutions are crepant and therefore
result in a smooth Calabi-Yau variety. 
We obtain these resolutions by triangulating the cone $\sigma$, as shown  in
Figure~\ref{torvar}(Right).
This is equivalent to constructing a lattice triangulation of the trapezoid:
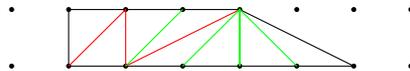
\begin{center}\begin{tikzpicture}[scale=.75]
\foreach \x in {-1,0,1,2,3,4,5,6} {
  \foreach \y in {0,1} {
    \node at (\x,\y) [circle,fill=black,scale=.2] {};
  }
}
\draw (0,0) -- (0,1) -- (3,1) -- (5,0) -- (0,0);
\draw[red] (0,0) -- (1,1) -- (1,0) -- (3,1);
\draw[green] (1,0) -- (2,1);
\draw[green] (2,0) -- (3,1) -- (3,0) -- (3,1) -- (4,0);
\end{tikzpicture}\end{center}
Internal edges in the triangulation of the strip correspond to
two-dimensional cones in the toric fan of the resolved threefold; they
describe the irreducible components of the exceptional curve. 
Absence of lattice points in the interior of the cone signifies that
the resolution does not contain compact $4$-cycles. Its second homology is
thus generated by the components of the exceptional curve.
Each prime component of the exceptional set is a smooth rational curve.

We shall consider all possible crepant resolutions of $C_{m,n}$, which
correspond to all maximal lattice triangulations of the strip
(i.e.\ triangulations in which each triangle has area $\frac12$). We
shall abuse notation to use $C_{m,n}$ to refer to the strip as well as
to the variety which it defines. Let us first collect some combinatorial
properties of these triangulations.
\begin{proposition}\label{count}\mbox{}\\[-\baselineskip]
\begin{enumerate}
\item Each triangulation of the polygon $C_{m,n}$ has $N_{F}=m+n$ triangles
and
$N_E= m+n-1$ interior edges.
\label{edg}
\item There are $N_{\triangle}=\binom{m+n}{m}$
triangulations of $C_{m,n}$.
\label{trgcnt}
\item The Euler characteristic of any crepant resolution of
$C_{m,n}$ is $m+n$.
\label{euler}
\end{enumerate}
\end{proposition}
\begin{proof}
The area of each regular triangle in a tesselation of the polygon
is $1/2$, as mentioned
above. The area of the trapezoid is $(m+n)/2$. Hence the number of triangles
in each triangulation is $N_F=m+n$. 

Since every interior edge of a triangulation emanates 
from a point in the upper row (also ends
on a point in the lower), it suffices to count the number of lines emanating
from the points in the upper row. Considering a triangulation,
let $N_i$  denote the number of triangles
containing the point $(i,1)$ in the upper row; $0\leq
i\leq m$. These triangles have a totality of
$\binom{N_i}{2}$ pairwise intersections
of which $\binom{N_i-1}{2}$ intersections are at the point alone,
while the rest, $\binom{N_i-1}{1}=N_i-1$ intersections, are along an interior 
edge, containing the point (cf. Figure~\ref{fig:intxn}). 
\begin{figure}[h]
\begin{pspicture}(-1,-1)(2,2)
\pspolygon(0,0)(3,0)(2,1)(0,1)
\pspolygon*(1,1)(0,0)(1,0)
\pspolygon*(1,1)(2,1)(2,0)
\psline(2,1)(2,0)
\rput(1,1.3){$\scriptstyle (1,1)$}
\end{pspicture}
\caption{The point $(1,1)$ is contained in four triangles. A pair of 
neighboring triangles intersect at an interior edge containing the point,
while non-neighbouring triangles intersect at the point only.}
\label{fig:intxn}
\end{figure}
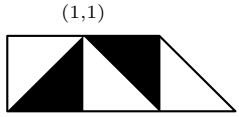
Hence the total number of interior edges
is $N_E=\sum_{i=0}^{m} (N_i-1)=\sum_{i=0}^{m}N_i - (m+1)$. On the other hand,
since the point $(i,1)$ is shared by $N_i$ triangles and there is a single
triangle containing two consecutive
points in the upper row, summing $N_i$ over all the points in the upper row
counts the number of triangles with $m$ triangles counted twice, ergo
$\sum_{i=0}^{m}N_i = N_F+m$. From these two expressions and the expression
for $N_F$ obtained above, we have $N_E = m+n-1$. This proves
statement~\ref{edg}.

To count the number of triangulations, let us note that all of the $N_E$
interior edges starts from one of the $m$ points in the upper row, which can
happen in $\binom{N_E}{m}$ ways. Also, all of these lines end on one of the 
$n$ points in the bottom row, which can happen in $\binom{N_E}{n}$ ways.
Adding, we have the number of triangulations $N_{\triangle} =
\binom{m+n-1}{m}+\binom{m+n-1}{n}=\binom{m+n}{m}$, proving
statement~(\ref{trgcnt}). 

Finally, for any crepant resolution the Euler characteristic is
$\chi(X) = h^0(X;\mathbb{Z}) + h^2(X;\mathbb{Z})$ 
in the absence of higher dimensional homology cycles. Moreover, for the cases
at hand, the two-cycles are given by the interior edges, so that
$ h^2(X;\mathbb{Z})=N_E$, while $ h^0(X;\mathbb{Z})=1$. 
Statement~\ref{euler} follows.

{\bf Aliter}:
We can count the number $N_F$ in another way by observing that
each triangle in a triangulation has a unique horizontal side,
which is either at the \emph{top} or at the \emph{bottom} of the strip,
corresponding to vertical coordinate $1$ or $0$, respectively. 
We shall refer to this side as the \emph{base of the triangle}. Since there
are $m$ segments on the top line and $n$ on the bottom, each of which 
is the base of one and only one triangle, the number of
triangles in a triangulation is $N_F=m+n$. 
\end{proof}
\subsection{Enumerating triangulations}
In the following we require a means to enumerating triangulations and 
labelling its triangles and edges.
There is a natural ordering
of triangles in a triangulation ``from left to right''.
We start with the unique
triangle $t_1$ having the line $(0,0)$--$(0,1)$ as its side
and move towards the right across the unique other non-horizontal edge to
arrive at the next triangle $t_2$. Continuing and labelling triangles on
the way seriatim, we finally   
arrive at the unique triangle $t_{m+n}$ which has the line
$(m,1)$--$(n,0)$ as its side. 
From the expression of $N_{\triangle}$ obtained above it is clear that 
specifying the
$m$ triangles based on the top line, or, alternatively, the $n$ triangles
based on the bottom line, fixes a triangulation. However, since we have
assumed $m\leq n$, the first choice is more economic and we shall adhere to
it. Hence, we denote each triangulation of
$C_{m,n}$ by a subset $T \subset \{1, 2, \dotsc, N_F\}$ with length
$\abs{T} = m$, where the base of the triangle $t_k$, $k \in T$ is
at the upper line
of the strip and $\{\oton[m+n]t\}$ denotes the set of all triangles.
These are illustrated in the following
\begin{example}
Let $m=2$ and $n=4$. Here are some  of the triangulations of $C_{2,4}$
given by subsets of length $2$ of $\bigl\{1,\dotsc,6\bigr\}$.
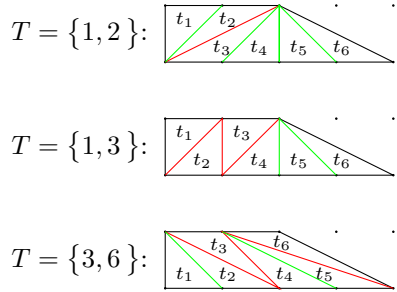
\begin{center}\begin{tikzpicture}[scale=.75]
\foreach \x in {2,3,4,5,6} {
  \foreach \y in {-1,0,1,2,3,4} {
    \node at (\x,\y) [circle,fill=black,scale=.1] {};
  }
}
\node at (.5, 3.5) { $T=\bigl\{1,2\,\bigr\}$: };
\draw (6,3) -- (4,4) -- (2,4) -- (2,3) -- (6,3);
\draw[red] (2,3) -- (4,4);
\draw[green] (2,3) -- (3,4);
\draw[green] (3,3) -- (4,4) -- (4,3) -- (4,4) -- (5,3);
\node at (2.35, 3.75) { $\scriptstyle t_1$ };
\node at (3.1, 3.75) { $\scriptstyle t_2$ };
\node at (3, 3.25) { $\scriptstyle t_3$ };
\node at (3.65, 3.25) { $\scriptstyle t_4$ };
\node at (4.35, 3.25) { $\scriptstyle t_5$ };
\node at (5.1, 3.25) { $\scriptstyle t_6$ };

\node at (.5, 1.5) { $T=\bigl\{1,3\,\bigr\}$: };
\draw (6,1) -- (4,2) -- (2,2) -- (2,1) -- (6,1);
\draw[red] (2,1) -- (3,2) -- (3,1) -- (4,2);
\draw[green] (4,1) -- (4,2) -- (5,1);
\node at (2.35, 1.75) { $\scriptstyle t_1$ };
\node at (2.65, 1.25) { $\scriptstyle t_2$ };
\node at (3.35, 1.75) { $\scriptstyle t_3$ };
\node at (3.65, 1.25) { $\scriptstyle t_4$ };
\node at (4.35, 1.25) { $\scriptstyle t_5$ };
\node at (5.1, 1.25) { $\scriptstyle t_6$ };

\node at (.5, -.5) { $T=\bigl\{3,6\,\bigr\}$: };
\draw (6,-1) -- (4,0) -- (2,0) -- (2,-1) -- (6,-1);
\draw[red] (2,0) -- (4,-1) -- (3,0) -- (6,-1);
\draw[green] (2,0) -- (3,-1);
\draw[green] (3,0) -- (5,-1);
\node at (2.35, -.75) { $\scriptstyle t_1$ };
\node at (3.10, -.75) { $\scriptstyle t_2$ };
\node at (2.95, -.25) { $\scriptstyle t_3$ };
\node at (4.1, -.75) { $\scriptstyle t_4 $};
\node at (4.8, -.75) { $\scriptstyle t_5$ };
\node at (4.05, -.2) { $\scriptstyle t_6 $};
\end{tikzpicture}
\end{center}
\end{example}
Interior edges are labelled using intersection of adjacent
triangles. We define the $i$-th edge
$e_i$ as
\[ e_i \ce t_i \cap t_{i+1} \text{ , $i=1, \dotsc, N_E$.} \]
In a given triangulation $T \subset \{1,\dotsc, N_F\}$, there are two
possibilities for each edge $e_i$, namely, it is either the intersection of 
two triangles $t_i$ and $t_{i+1}$ both having bases on the same
horizontal line (top or bottom) of the strip, or they have bases
on different lines. In the former case
either $i,i+1 \in T$ or $i,i+1\not\in T$, 
 we say that $e_i$ is of type
``$+$'' and colour the edge green in the toric diagram. These correspond to
$\mathcal{O}(-2,0)$ curves.
In the latter case either $i \in T, i+1 \not\in T$ or $i \not\in T, i+1 \in T$,
$e_i$ is said to be of type ``$-$'' and we depict it in red. These correspond
to $\mathcal{O}(-1,-1)$ curves.  
We let $\tau(e_i) = \pm 1$
according to whether $e_i$ is of type ``$+$'' or ``$-$''.
This furnishes a canonical scheme for ordering and characterizing the edges,
which correspond to bases of the second cohomology group of crepant
resolutions. 
\subsection{Computing triangulations}
In working with triangulations implementation of the above scheme in
computer programs is useful. Let us briefly discuss some aspects.
The triangulation was carried out using the software TOPCOM
\cite{TOPCOM}. The function \texttt{points2allfinetriangs},
triangulates a strip using triangles of equal, minimal area 
producing a list of all possible triangulations.

In TOPCOM, points in a point set are given in
homogeneous coordinates, so for our purposes the vertex $(i,j)$
corresponds to the point \texttt{[$i$,$j$,1]}.  We label the $m+n+2$
vertices sequentially, assigning the range $0, \dotsc, m$ to the
vertices $v_0 \ce\,$\texttt{[0,0,1]}, $v_1 \ce\,$\texttt{[1,0,1]}, \ldots, $v_m \ce\,$\texttt{[$m$,0,1]},
and the range $m+1, \dotsc, m+n+1$ to $v_{m+1} \ce\,$\texttt{[0,1,1]}, $v_{m+2} \ce\,$\texttt{[1,1,1]}
\ldots, $v_{m+n+1} \ce\,$\texttt{[$n$,1,1]}. The output of
TOPCOM consists of lists of triplets $(v_a,v_b,v_c)$ of vertices
giving the triangulation of the strip. 
The internal edges in a triangulation are extracted from this list.

Their types are determined as follows.
The natural ordering ``from left to right'' of the non-horizontal
edges is precisely the lexicographic ordering of either the top or
the bottom vertices $(i,j)$. When the edges are ordered in this fashion,
the $k$-{th} edge, corresponding to the 
vertex $(i_k, j_k)$, is of type ``$+$''
if $j_{k-1} = j_k = j_{k+1}$ and $i_{k-1} + 1 = i_k = i_{k+1} - 1$. 
Otherwise it is of type ``$-$'' 
\footnote{We are grateful to Jesus Martinez-Garcia for writing
the program.}.

From this data we can construct the partition function of any particular
resolution of $C_{m,n}$ corresponding to a specific triangulation.


\section{Curve counting on singular varieties}
For any complex threefold $(X,\mathcal{O}_X)$, the Hilbert scheme of
ideal sheaves $\mathcal{I} \subset \mathcal{O}_X$ with fixed Euler
characteristic $\chi(\mathcal{I}) = k$ and support
$[\operatorname{supp} (\mathcal{I})] = \beta \in H_2(X;\mathbb{Z})$,
written $I_k(X,\beta)$, has a perfect obstruction theory of virtual
dimension $\int_\beta c_1(T_X) = -K_X(\beta)$, see \cite{DT}. When
$K_X=0$, the numbers
\[ N_{k,\beta}(X) \ce \int_{[I_k(X,\beta)]^\text{vir}} 1 \]
are the \emph{Donaldson--Thomas (DT) invariants} of $X$. Let 
$Q=(Q_1,\dotsc,Q_h)$, $h=\text{dim} H_2(X,\Z)$, 
be a set of symbols corresponding to generators of $H_2(X;\mathbb{Z})$. 
The DT invariants
are collected into the \emph{Donaldson--Thomas partition function}
\[ Z(X;q,Q) \ce \sum_{k=0}^\infty \ \sum_{\beta \in H_2(X;\mathbb{Z})}
   N_{k,\beta}(X) \; q^k \; Q^\beta, \]
where $Q^{\beta} = Q_1^{\beta_1} \dotsm Q_h^{\beta_h}$. We single out the degree-$0$ contributions,
\[ Z_0(X;q) \ce \sum_{k=0}^\infty N_{k,0}(X) \; q^k \text{ ,} \]
and we define the \emph{reduced DT partition function} as
\[ Z'(X;q,Q) \ce Z(X;q,Q) / Z_0(X;q). \]

For any smooth, toric threefold $X$, we have $K_X(0) = 0$ and so we
can define the degree-$0$ partition function $Z_0(X;q)$. It is known
\cite{MNOP1} that
\[ Z_0(X; -q) = M(1,q)^{\textstyle\int_X c_3(T_X \otimes K_X)} \text{ ,} \]
and in particular if $X$ is Calabi-Yau, then
\[ Z_0(X; -q) = M(1,q)^{\chi(X)}, \]
where $\chi(X)$ denotes the Euler characteristic of $X$ and 
\begin{equation}
\label{mac} 
M(x,q) \ce \prod_{k=1}^\infty \frac{1}{(1-xq^k)^k} = \exp
\sum_{\vphantom{j}i=1}^\infty \sum_{j=1}^\infty \frac{i}{j}~x^j~q^{ij}. 
\end{equation} 
denotes the (generalised) MacMahon function. The
nexus between the partition function and the MacMahon function originates 
from the fact that the MacMahon function counts box partitions, and
degree-$0$ toric ideal sheaves are given by monomial ideals, which can
indeed be arranged like ``boxes stacked into a corner''.
\subsection{DT invariants of generalized conifolds}
If $X$ is a crepant resolution of $C_{m,n}$, then it is a smooth, toric
Calabi-Yau threefold. The DT partition function can be computed
combinatorially by the topological vertex method (see \cite{LLLZ,IKP}).
We shall always take the curves corresponding to the interior
edges $e_i$ as our preferred basis for $H_2(X;\mathbb{Z})$, that is,
\[ \beta = \sum_{i=1}^{N_E} \beta_i [e_i] \in H_2(X;\mathbb{Z}), \]
where $N_E=m+n-1$, by Proposition~\ref{count}.
Furthermore, we have $\chi(X) = m+n$.

We need to establish some terminology to describe $Z'(X;q,Q)$. A
set $P = \{i,i+1, \dotsc, j\}$ is called 
an \emph{edge path} if $1\leq i \leq j \leq
N_E$. It is to be thought of as a sequence of consecutive interior edges
of the triangulation $T$ of $C_{m,n}$ corresponding to the resolution
$X$. An edge path $P$ has \emph{length} $\abs{P} \ce j-i+1$ and connects the
triangles $t_i$ and $t_{j+1}$.  In a triangulation there
are $m+n-1$ edge paths of length $1$, $m+n-2$ of length $2$, and so
forth, and $1$ of length $m+n-1$, so in total there are
$\binom{m+n}{2}$ edge paths.  An edge path is literally a path along
the compact edges of the dual tropical curve of the triangulation $T$.

If $P = \{i,i+1,\dotsc,j\}$ is an edge path, we write $Q_P = Q_{ij} =
Q_i \dotsm Q_j$, so for example $Q_{22} = Q_2$ and $Q_{35} = Q_3Q_4Q_5$.
We define
\[ f(P, q, Q) = M(Q_P, q)^{\tau(e_i) \tau(e_{i+1}) \dotsm \tau(e_j)}. \]
Thus, $f(P,q,Q)$ is either the MacMahon function or its reciprocal,
depending on whether $P$ contains an even or an odd number of edges of
type ``$-$''. The whole partition function of $X$ is  the
product of such terms over all edge paths, that is,
\begin{equation} 
\label{prod}
Z'(X; -q, Q) = \prod_{P} f(P,q,Q) 
= \prod\limits_{\substack{\{i,j|\\1\leq i\leq j\leq N_E\}}}
\prod\limits_{k=1}^{\infty} 
\left( 1 - Q_{ij} q^k\right)^{-k\tau(e_i)\cdots\tau(e_j)}
\end{equation} 
Since this partition function is determined entirely by the
triangulation, i.e.\ by the subset $T \subset \{1, 2, \dotsc, N_F\}$,
$\abs{T} = m$, alluded to above,
we write $Z'_T(C_{m,n}; q,Q^T)$ for the partition
function, where we  abbreviate
$Q^T = (\oton[N_E]{Q^T})$. We now consider the
collection of all possible triangulations of $C_{m,n}$.
\begin{definition}\label{def.tot}
We define the \emph{total partition function}:
\[ Z'_\text{tot}(C_{m,n}; -q, Q) \ce \prod_{\substack{T \subset
\{1,2,\dotsc,N_F\} \\ \abs{T} = m}} Z'_T(C_{m,n}; -q, Q) \text{.} \]
\end{definition}
Let us consider the following ad-hoc definition.
\begin{definition}\label{def.hom}
A partition function $Z(q,Q)$ of variables $Q = (Q_1, Q_2, \dotsc)$ is
called \emph{homogeneous} if
\[ Z(q, Q) = \left( \prod M\bigl(\textstyle\prod_{i \in A} Q_i, q\bigr) \right)^d, \]
where the first product is over an arbitrary finite collection of index set $A \subset \{1, 2, \dotsc\}$.
The exponent $d$ is called the \emph{degree} of $Z$.
\end{definition}
\begin{example}
\label{ex:predef}
Let us consider $m=1$, $n=1$, in which case 
the strip is a single square admitting two triangulations, namely,
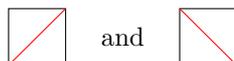
\begin{center}\begin{tikzpicture}[scale=.75]
\draw (2,3) -- (3,3) -- (3,4) -- (2,4) -- (2,3);
\draw[red] (2,3) -- (3,4);

\node at (4,3.5) { and };

\draw (5,3) -- (6,3) -- (6,4) -- (5,4) -- (5,3);
\draw[red] (5,4) -- (6,3);

\node at (7,3.5) { };
\end{tikzpicture}
\end{center}
yielding the partition function $Z'_\text{tot}(C_{1,1}; -q, Q) =
M(Q_1,q)^{-2}$, which is homogeneous with degree $-2$.
Triangulations on smaller strips can be extended to triangulations 
of bigger strips. To illustrate this 
let us consider the following two ways to pass 
from a triangulation of $C_{m,n}$ to a
triangulation of $C_{m,n+1}$. In the first case, the right-most edge
of $C_{m,n}$ turns into an internal edge of $C_{m,n+1}$ of ``$+$''
type, as
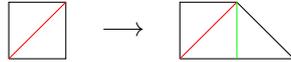
\begin{center}\begin{tikzpicture}[scale=.75]
\draw (2,3) -- (3,3) -- (3,4) -- (2,4) -- (2,3);
\draw[red] (2,3) -- (3,4);

\node at (4,3.5) { $\longrightarrow$ };

\draw (5,3) -- (7,3) -- (6,4) -- (5,4) -- (5,3);
\draw[red] (5,3) -- (6,4);
\draw[green] (6,3) -- (6,4);
\end{tikzpicture}
\end{center}
The exponent of $M(Q_{1,m+n-1},q)$ coming from this triangulation of
$C_{m,n}$ is the same as the exponent of $M(Q_{1,m+n},q)$ for the
corresponding triangulation of $C_{m,n+1}$. Hence, there is a
correspondence between such kinds of triangulations of the two strips,
maintaining equality of exponents of the MacMahon factors. In the
second case the rightmost edge of $C_{m,n}$ turns into an internal
edge of $C_{m,n+1}$ of ``$-$'' type, as
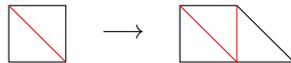
\begin{center}\begin{tikzpicture}[scale=.75]
\draw (2,3) -- (3,3) -- (3,4) -- (2,4) -- (2,3);
\draw[red] (2,4) -- (3,3);

\node at (4,3.5) { $\longrightarrow$ };

\draw (5,3) -- (7,3) -- (6,4) -- (5,4) -- (5,3);
\draw[red] (5,4) -- (6,3) -- (6,4);
\end{tikzpicture} 
\end{center}
Now in the triangulation on the right-hand side we 
have MacMahon factors as $M(Q_1,q)^{-1}$ and $M(Q_1Q_2,q)^{+1}$, which
appears to give rise to different exponents. However, since every
parallelogram has two diagonals, there is another triangulation
obtained by flopping the diagonal on rightmost parallelogram of the
previous figure, and we obtain an extra triangulation of $C_{m,n+1}$
(this one not coming from a triangulation of $C_{m,n}$) as
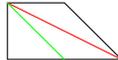
\begin{center}\begin{tikzpicture}[scale=.75]
\draw (2,3) -- (4,3) -- (3,4) -- (2,4) -- (2,3);
\draw[red] (2,4) -- (4,3);
\draw[green] (2,4) -- (3,3);
\end{tikzpicture}
\end{center}
that contributes factors of 
$M(Q_1,q)^{+1}$ and $M(Q_1Q_2,q)^{-1}$, cancelling out the seemingly 
unbalanced contributions from the previous one. 
\end{example}
In general, we have
\begin{proposition}\label{prop.homo}
For $0 < m \leq n$, $Z'_\text{\rm tot}(C_{m,n}; -q, Q)$ is homogeneous of degree $d$,
where
\begin{equation}
\label{eq.a}
d = \frac{(m^2-m+n^2-n-2mn)(m+n-2)!}{m!n!},
\end{equation}
namely,
\[ Z'_\text{\rm tot}(C_{m,n}; -q, Q) 
= \prod_{1\leq i \leq j \leq m+n-1} M(Q_{ij},q)^d,\]
\end{proposition}
\begin{proof}
We first present a purely combinatorial proof.
The proposition consists of two separate parts, and so does the
proof. The first statement is that each MacMahon factor $M(Q_{ij},q)$
appears with the \emph{same} power in the total partition
function.

We have to show that each MacMahon factor $M(Q_{ij},q)$ appears with
the \emph{same} power in the total partition function and compute the value of
this exponent. The problem is entirely combinatorial. In terms of
finite sets, it takes the following form: Let us simply write $N$ for
the finite set $\{1,2,\dotsc,N\}$. For any subset $T \subset N$ and
any fixed, ordered subset $S = \{\oton[k]s\} \subset N$, we define the
\emph{characteristic sequence}
\[ \chi_T(S) \ce \bigl(\chi_T(s_1), \dotsc, \chi_T(s_k)\bigr), \]
where $\chi_T \colon N \to \{0,1\} = \mathbb{Z}/2\mathbb{Z}$ is the
characteristic function of $T$. (It will be opportune to think of the
two-element set as the additive group of order $2$.)

In our application, we shall take $S$ to be a ``contiguous'' subset of
the form $\{i, i+1, \dotsc, j\}$ corresponding to some edge path. For such a subset,
we define the \emph{difference sequence} as
\[ \Delta_T(S) \ce \bigl(\chi_T(s_1) - \chi_T(s_2), \chi_T(s_2) - \chi_T(s_3), \dotsc,
   \chi_T(s_{k-1}) - \chi_T(s_k)\bigr) \text{ ,} \]
and we define the \emph{$T$-signature} of $S$ as
\[ \sigma_T(S) \ce \prod_{b \in \Delta_T(S)} (-1)^b \in \{+1,-1\} \text{ .} \]
(Since we are only interested in the $T$-signature, we may consider the elements
of $\Delta_T(S)$ to take values in $\mathbb{Z}/2\mathbb{Z}$ and identify $+1$ and $-1$.)
Finally, the exponent of $M(Q_{ij}, q)$ in the total partition function of $C_{m,n}$
is the \emph{$m$-signature} of the set $S = \{i,i+1,\dotsc,j\}$, defined as
\[ \sigma(S) = \sum_{T \subset N \;:\; \abs{T} = m} \sigma_T(S) \text{ ,} \]
where $N = m+n$.

So much for the setup. The first observation is that any action $\pi
\in \Sigma_N$ that preserves the contiguous ordering of the elements of $S$
does not alter the value of the total signature: $\sigma(\pi S) =
\sigma(S)$. Therefore, we may assume without loss of generality that
$S$ is $\{1,2,\dotsc,k\}$.

Next, any subset $T \subset N$ with $\abs{T}=m$ is of the form $T = U
\sqcup T'$, where $U \subset \{1,2,\dotsc,k\}$ with $\abs{U}=i$ and
$T' \subset \{k+1,k+2,\dotsc,N\}$ with $\abs{T'}=m-i$ for
$i=0,1,\dotsc,k$. Now observe that all we need to compute the
$m$-signature is $\Delta_U(S)$, or rather $\sigma_U(S) = \sigma_T(S)$.
Since there are $\binom{N}{m}$ subsets in total,
we have
\[ \sigma(S) = \abs{\{ T : \sigma_T(S) = +1 \bigr\}} - \abs{\{ T : \sigma_T(S) = -1 \bigr\}}
    = \binom{N}{m} - 2 \abs{\{ T : \sigma_T(S) = -1 \bigr\}} \text{ .} \]
The combinatorics of this are easily determined: Subsets $T = U \sqcup T'$ for which $\sigma_U(S) = -1$
are those for which $\Delta_U(S)$ has an odd number of $1$s, and there are $2\binom{k-2}{i-1}$ of those,
where $i = \abs{U}$. Summing over all $i$ we find
\begin{equation*} 
\sigma(S) = \binom{N}{m} - 4 \sum_{i=1}^{k-1} \binom{k-2}{i-1} 
\binom{N-k}{m-i}.
\end{equation*} 
The last factor accounts for all the possible subsets $T'$. 
The sum evaluates to $\binom{N-2}{m-1}$ leading to
\[ \sigma(S) = \binom{N}{m} - 4 \binom{N-2}{m-1} =
\frac{(N^2-N+4m^2-4mN)(N-2)!}{m!(N-m)!}. \]
This is true for any contiguous set $S=\{i,i+1,\dotsc,j\}$, and the result follows by substituting $N=m+n$.

The second statement is the value of the exponent. Since the exponent
is the same for each factor $M(Q_{ij},q)$ by the first part, we may
compute it by just computing the exponent of $M(Q_1,q)$, i.e.\ the
factor corresponding to the edge path $\{1\}$. Each triangulation $T$
contributes either an exponent $+1$ or $-1$. The exponent is $+1$ if
$1,2 \in T$ or $1,2 \not\in T$, and it is $-1$ if $1 \in T$, $2
\not\in T$ or if $1 \not\in T$, $2 \in T$. The number of $+1$s is thus
the sum of the number of triangulations of $C_{m-2,n}$ and $C_{m,n-2}$,
and the number of $-1$s is twice the number of triangulations of $C_{m-1,n-1}$.

{\bf Aliter:} We present another proof using probabilities.
As discussed before, the interior edges $e_i$,  
$1,2,\cdots , N_E$ are  numbered from left to  right in a unique
fashion. An edge path $P={i,\cdots ,j}$ is then a sequence of
interior edges from $e_i$ to $e_j$, both inclusive,
$1\leq i\leq j\leq N_E$ connecting two
triangles $t_i$ and $t_{j+1}$. This is also
illustrated in Figure~\ref{nspic}. 
An edge path contributes a MacMahon
factor with a positive exponent 
to the partition function if it connects two triangles which are
either both based on the top line or both based on the bottom line. If it
connects two triangles based on different lines, then the contribution to the
partition function comes with a negative exponent. 

For a triangulation,
given an edge path $Q_{ij}$, the probability that the triangle $t_i$ has 
its base on the top line is $m/(m+n)$, as there are $m$ triangles with bases
on the top line in any triangulation and there are $m+n$ triangles in total. 
Then the probability that the triangle $t_{j+1}$ has also its base 
on the top line is
$(m-1)/(m+n-1)$. The probability that the edge path connects two triangles
both having bases on the top line is thus $p_t=m(m-1)/(m+n)(m+n-1)$. Similarly,
the probability that an edge path connects two triangles both having bases 
on the
bottom line is $p_b = n(n-1)/(m+n)(m+n-1)$. The probability that an edge path
connects triangles having bases on different lines is then $1-p_t-p_b$. Hence
considering all the $N_{\triangle}$ triangulations, the
contribution to the partition function comes with the exponent
\begin{equation} 
\label{eq:expo}
\begin{split}
d &= \left(p_t + p_b - (1-p_t-p_b)\right)N_{\triangle}\\
  &= \left(\frac{2m(m-1)}{(m+n)(m+n-1)}+\frac{2n(n-1)}{(m+n)(m+n-1)}-1\right)
\frac{(m+n)!}{m!n!}\\
&= 2\binom{m+n-2}{m-2} + 2\binom{m+n-2}{n-2}-\binom{m+n}{n}\\
&= \frac{(m^2-m+n^2-n-2mn)(m+n-2)!}{m!n!}.
\end{split}
\end{equation} 
While the integrality of the  exponent $d$ is obvious from its definition,
we made it conspicuous by writing it 
as a combination of binomial coefficients in the third line.
\end{proof}
\begin{remark}
The case $n>m=0$ is excluded from the first proof of the proposition,
since $C_{0,n}$
only admits one unique triangulation, and all interior edges are of
type ``$+$''. Writing $X$ for the resolution, we have
\[ Z'(X; -q, Q) = \!\!\!\!\!\!\! \prod_{1\leq i \leq j \leq n-1} \!\!\!\!\!\!\! M(Q_{ij},q)
   \eqand Z(X; -q, Q) = M(1,q)^n \; Z'(X; -q, Q) \text{ .} \]
We have indeed $d=1$ in equation~\ref{eq.a} whenever $m=0$.

The second proof, on the other hand,
only excludes the case $m=0$, $n=1$, for not having
any interior edge. It is more general in this sense.
\end{remark}
\begin{definition}\label{cc}
A partition function for a Calabi-Yau manifold $Y$
is of {\em curve-counting type} if it can be expressed in terms of
the Donaldson--Thomas partition function up to a factor depending
only on the Euler characteristic of $Y$.
\end{definition}
 We have thus proved:

\begin{theorem}\label{th.homo}
Let $X$ be a toric singular Calabi-Yau threefold without contractible curves
or compact $4$-cycles.
Let $Z(Y;q,Q)$ be any partition
function of curve-counting type. Then the total partition function
for $X$ is given by 
\[ Z_\text{\rm tot}(X;q,Q) \ce \prod_{\substack{Y\\Y \to X}} 
Z(Y;q,Q), \]
where the product ranges over all crepant resolutions of $X$, is homogeneous,
and its degree is given by Proposition~\ref{prop.homo}.
\end{theorem}

For a general singular toric Calabi-Yau threefold $X$ without compact
$4$-cycles, we can use this theorem to factor the partition function
into homogeneous factors. The toric diagram $\Delta$ of $X$ is a strip
of shape $C_{m,n}$ with an arbitrary number of internal edges filled
in, for example,
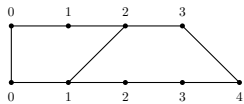
\begin{center}\begin{tikzpicture}[scale=.75]
\draw  (0,0) -- (0,1)  -- (3,1) -- (4,0) -- (0,0);
\draw (1,0) -- (2,1);
\foreach \x in {0,1,2,3} {
  \foreach \y in {0,1} {
    \node at (\x,\y) [circle,fill=black,scale=.2] {};
  }
}
\node at (4, 0) [circle,fill=black,scale=.2] {};
\foreach \x in {0,1,2,3,4} { \node at (\x, -.25) [scale=.5] {\x}; }
\foreach \x in {0,1,2,3}   { \node at (\x, 1.25) [scale=.5] {\x}; }
\end{tikzpicture}
\end{center}
Let us partition the integers $m, n$ according to the already 
filled-in interior edges, that is,
\[ (m,n) = \sum_{k=1}^{P} (m_k,n_k) = (m_1 + m_2 + 
\dotsb + m_P, n_1 + n_2 + \dotsb + n_P). \]
In the example above, we have $(m,n) = (3,4)$, and the single interior
edge corresponds to the partition $(3,4) = (2+1, 1+3)$. It is clear
that the number of maximal triangulations of this shape is
\[ \prod_{k=1}^P \binom{m_k+n_k}{n_k} \text{ ,} \]
where each factor counts the number of triangulations of the embedded
subdiagram $C_{m_k,n_k} \ec C_k$. If we restrict our attention to some
fixed subdiagram $C_k$, then entire collection of triangulations of
$\Delta$ contains many triangulations with the same restriction to
$C_k$. It is clear that for any fixed triangulation of $C_k$, there
are $b_k$ triangulations of $\Delta$ that restrict to the given triangulation,
where
\[ b_k = \prod_{j \neq k} \binom{m_j+n_j}{n_j}. \]
We extend Definition \ref{def.tot} in a straightforward manner to
\begin{definition}
If $X$ is a singular Calabi-Yau threefold without compact $4$-cycles such that
the convex hull of its toric diagram is $C_{m,n}$ (that is, there
exists a birational map $X \to C_{m,n}$), we define the \emph{total
partition function} to be
\[ Z'_\text{tot}(X; -q,Q) \ce \prod_T Z'_T(C_{m,n}, -q, Q) \text{ .} \]
Here the term in the product of the right-hand side is the same as in
Definition \ref{def.tot}, except that the product is taken only over
those triangulations $T$ which correspond to resolutions of $X$.
\end{definition}
Now Theorem \ref{th.homo} implies the following:
\begin{corollary}\label{contrac}
If $X$ is a singular Calabi-Yau threefold  without compact $4$-cycles and
$(m,n)$, $P$ and $b_k$ are as above, then the total partition function of $X$
factors as follows:
\[ Z'_\text{\rm tot}(X; -q, Q) = Z''(-q,Q) 
\prod_{k=1}^P Z'_\text{\rm tot}(C_{m_k,n_k}; -q, Q)^{b_k}.\]
The factors in the product on the right are homogeneous as per Theorem~\ref{th.homo},
and the function $Z''$ only contains factors $M(Q_{ij},q)$ for which
the edge path corresponding to $Q_{ij}$ crosses one of the interior
edges of the toric diagram of $X$.
\end{corollary}
\begin{example}
In the above example with $(m,n) = (3,4) = (2+1,1+3)$, 
the two homogeneous factors are $Z'_\text{tot}(C_{1,2};-q,Q)^3$
and $Z'_\text{tot}(C_{3,1};-q,Q)^2$, and the inhomogeneous
factor contains only terms $M(Q_{ij},q)$ with $i \leq 3 \leq j$,
because the third edge is already fixed in the diagram.
\end{example}

\subsection{BPS counting and relation to black holes}
Here is one application to BPS state counting. The \emph{topological
string partition function} of $X$ is
\[ Z_\text{top}(X; q, Q) = M(1,q)^{\chi(X)/2} Z'(X;-q,Q) \text{ ,} \]
so it is a partition function of curve-counting type.

\begin{corollary}\label{BH}
Writing $X_T$ for the resolution of $C_{m,n}$ corresponding to the
triangulation $T$, we have
\[ \prod_{T} Z_\text{top}(X_T;q,Q) = M(1,q)^{\binom{m+n}{m}\frac{m+n}{2}}
   \prod_{1 \leq i \leq j \leq m+n-1} M(Q_{ij},q)^{\frac{(m^2-m+n^2-n-2mn)(m+n-2)!}{m!n!}} \text{ .} \]
\end{corollary}
\begin{proof}
This follows immediately from the fact that $\chi(X_T) = m+n$ for all $T$ and
that there are $\binom{m+n}{m}$ triangulations.
\end{proof}
\section{Partition function via change of variables}
\label{sec:app}
\begin{figure}[h]
\psset{linewidth=.4pt}
\begin{pspicture}(-.5,-1)(14,2)
\pspolygon(0,0)(3,0)(2,1)(0,1)
\psline(1,1)(0,0)
\psline(1,1)(1,0)
\psline(1,1)(2,0)
\psline(2,1)(2,0)
\rput(.25,.75){\ns}
\rput(.75,.25){\nsp}
\rput(1.25,.25){\nsp}
\rput(1.75,.75){\ns}
\rput(2.25,.25){\nsp}
\pspolygon(5,0)(8,0)(7,1)(5,1)
\psline(5,1)(6,0)
\psline(6,1)(6,0)
\psline(6,1)(7,0)
\psline(7,1)(7,0)
\rput(5.25,.25){\nsp}
\rput(5.75,.75){\ns}
\rput(6.25,.25){\nsp}
\rput(6.75,.75){\ns}
\rput(7.25,.25){\nsp}
\pspolygon(10,0)(13,0)(12,1)(10,1)
\psline(10,1)(11,0)
\psline(10,1)(12,0)
\psline(10,1)(13,0)
\psline(11,1)(13,0)
\rput(10.25,.25){\nsp}
\rput(11,.25){\nsp}
\rput(11.85,.25){\nsp}
\rput(10.95,.85){\ns}
\rput(11.85,.85){\ns}
\end{pspicture}
\begin{pspicture}(-1,-2.5)(14,1)
\rput(.25,.50){\ns}
\rput(.75,.50){\nsp}
\psline(.25,.50)(.75,.50)
\rput(1.25,.5){$\scriptstyle ~Q_1(-)$}
\rput(.75,0){\nsp}
\rput(1.25,0){\nsp}
\psline(.75,0)(1.25,0)
\rput(1.75,0){$\scriptstyle ~Q_2(+)$}
\rput(1.25,-.50){\nsp}
\rput(1.75,-.50){\ns}
\psline(1.25,-.50)(1.75,-.50)
\rput(2.25,-.50){$\scriptstyle ~Q_3(-)$}
\rput(1.75,-1){\ns}
\rput(2.25,-1){\nsp}
\psline(1.75,-1)(2.25,-1)
\rput(2.75,-1){$\scriptstyle ~Q_4(-)$}
\rput(.25,-1.50){\ns}
\rput(1.25,-1.50){\nsp}
\psline(.25,-1.5)(1.25,-1.5)
\rput(1.85,-1.50){$\scriptstyle ~Q_{12}(-)$}
\rput(.25,-2){\ns}
\rput(1.75,-2){\ns}
\psline(.25,-2)(1.75,-2)
\rput(2.35,-2){$\scriptstyle ~Q_{13}(+)$}
\rput(5.25,.50){\nsp}
\rput(5.75,.50){\ns}
\psline(5.25,.50)(5.75,.50)
\rput(6.75,.5){$\scriptstyle ~R_1=1/Q_1(-)$}
\rput(5.75,0){\ns}
\rput(6.25,0){\nsp}
\psline(5.75,0)(6.25,0)
\rput(6.75,0){$\scriptstyle ~R_2(-)$}
\rput(6.25,-.50){\nsp}
\rput(6.75,-.50){\ns}
\psline(6.25,-.50)(6.75,-.50)
\rput(7.25,-.50){$\scriptstyle ~R_3(-)$}
\rput(6.75,-1){\ns}
\rput(7.25,-1){\nsp}
\psline(6.75,-1)(7.25,-1)
\rput(7.75,-1){$\scriptstyle ~R_4(-)$}
\rput(5.25,-1.50){\nsp}
\rput(6.25,-1.50){\nsp}
\psline(5.25,-1.5)(6.25,-1.5)
\rput(6.85,-1.50){$\scriptstyle ~R_{12}(-)$}
%
\rput(5.75,-2){\ns}
\rput(6.75,-2){\ns}
\psline(5.75,-2)(6.75,-2)
\rput(7.35,-2){$\scriptstyle ~R_{23}(+)$}
\rput(10.25,.50){\nsp}
\rput(10.75,.50){\nsp}
\psline(10.25,.50)(10.75,.50)
\rput(11.25,.5){$\scriptstyle ~S_1(+)$}
\rput(10.75,0){\nsp}
\rput(11.25,0){\nsp}
\psline(10.75,0)(11.25,0)
\rput(11.75,0){$\scriptstyle ~S_2(+)$}
\rput(11.25,-.50){\nsp}
\rput(11.75,-.50){\ns}
\psline(11.25,-.50)(11.75,-.50)
\rput(12.25,-.50){$\scriptstyle ~S_3(-)$}
\rput(11.75,-1){\ns}
\rput(12.25,-1){\ns}
\psline(11.75,-1)(12.25,-1)
\rput(12.75,-1){$\scriptstyle ~S_4(+)$}
\rput(10.25,-1.50){\nsp}
\rput(11.25,-1.50){\nsp}
\psline(10.25,-1.5)(11.25,-1.5)
\rput(11.85,-1.50){$\scriptstyle ~S_{12}(+)$}
\rput(10.25,-2){\nsp}
\rput(11.75,-2){\ns}
\psline(10.25,-2)(11.75,-2)
\rput(12.35,-2){$\scriptstyle ~S_{13}(-)$}
\end{pspicture}
\caption{NS-NS' brane arrangements corresponding to certain triangulations of
$C_{2,3}$. Filled circles denote NS-branes and white ones denote NS'-branes.
The path edges and their exponent in the partition function are
also indicated.}
\label{nspic}
\end{figure}
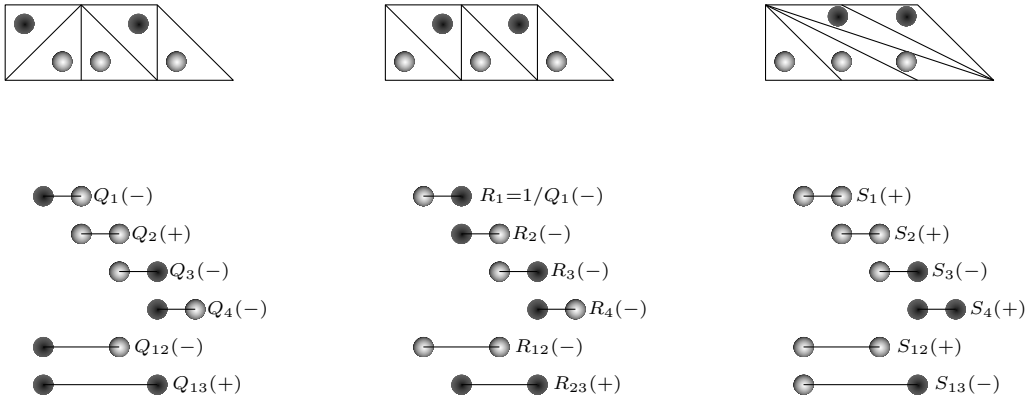
It has been mentioned earlier that 
the product of partition functions corresponding 
to different triangulations depend on the explicit isomorphism between 
homologies of crepant resolutions. For purposes of comparison let us 
briefly discuss the product of partition functions in the case when the 
map between the homologies of crepant  resolutions in different chambers 
in the moduli space is given by Seiberg duality \cite{Sz,Na,Y}.
We shall consider  the combinatorial aspects of the partition function in
terms of the dual type-IIA picture, given by a gauge theory 
of NS five-branes with D4-branes
stretched between them, interpreted as fractional branes. 
Depending on the spatial directions occupied by the 
NS branes in the target space, two types of branes,
referred to as NS and NS' branes, are considered.
The arrangement of the two types of 
NS-branes on a circle corresponds to the triangulations of the trapezoidal
strip \cite{uranga}. The field theory of such configuration of branes
is well-developed \cite{uranga,mp,unge}. We shall not discuss the field
theory here but focus only on certain combinatorial aspects of arrangement of
branes.

The T-dual type-IIA  theory on $C_{m,n}$
has $m$ NS branes and $n$ NS' branes. For any
triangulation of the trapezoid an NS brane corresponds to a regular triangle 
based on the top line and we denote it by a dark circle in Figure~\ref{nspic}. 
An NS' brane, on the
other hand, corresponds to a triangle based on the bottom line and will be
denoted by a white circle. The D4-branes stretched between these 
are denoted by a line, which also serves to designate
the relative separation between the NS-branes, given by the
period of the B-field. 

Considering an arrangement of NS-branes, a pair of branes linked by a line
corresponds to an edge path and contributes a factor to $Z'(X;-q,Q)$ in
\eqref{prod}. According to the combinatorial rule laid out earlier,
the index of the factor is positive if the branes are of the same
type, that is, the edge path connects 
either an NS-NS or an NS'-NS' pair and negative
otherwise. Indeed, a curve connecting two adjacent cones in the toric
diagram is $\mathcal{O}(-2)$ if the branes in the cones are of the 
same type and is $\mathcal{O}(-1,-1)$ otherwise. 
Thus, in particular, an NS and an
NS' branes are exchanged under a flop, as in $Q_1$ to $R_1$ in
Figure~\ref{nspic}.
For example, $Q_1$ in Figure~\ref{nspic} 
contributes
$\prod_{k=1}^{\infty}(1-q^kQ_1)^{-k}$ to the partition function, 
as it connects branes of
different types. On the other hand, $Q_{13} = Q_1Q_2Q_3$ contributes a factor
of $\prod_{k=1}^{\infty}(1-q^kQ_{13})^k$ as it connects branes 
of the same type. 

In this dual theory each triangulation of the strip corresponds 
to a ``phase" of the field theory described by a quiver gauge theory 
with a superpotential.
Different phases correspond to different paths to
approach the singularity from the asymptotic large-volume region. 

We consider products over all crepant resolutions, that is, phases, again
and we still assume all relevant resolutions to have isomorphic second 
homologies. 
However, instead
of formally identifying the elements of $H_2$,
we change coordinates to write each element of $H_2(X^s,\mathbb Z)$
in terms of a fixed basis $Q^o$.
We set
\[ Z^a_\text{new}(X;Q,\dotsc) \ce 
\prod_{t \in \mathcal{T}} Z_\text{old}(X^t; Q^t(Q^o), \dotsc). \]
\begin{example}\label{basicflop}
For the conifold we have two crepant resolutions
\begin{center}\begin{tikzpicture}[scale=.75]
\draw (2,3) -- (3,3) -- (3,4) -- (2,4) -- (2,3);
\draw [blue] (2,3) -- (3,4);

\node at (4,3.5) { and };
\node at (1.5,3.5) { Q };
\node at (5.5,3.5) { R } ;

\draw (6,3) -- (7,3) -- (7,4) -- (6,4) -- (6,3);
\draw (6,4) -- (7,3);

\node at (7,3.5) { };
\end{tikzpicture}\text{,}
\end{center}
with respective basis for the second homologies denoted $\{Q_1\}$ 
and $\{R_1\}$. Since they are related by a flop,
the change of coordinates reads $R_1=Q_1^{-1}$. We obtain the partition
function \cite{Sz}, 
\[ Z^a_\text{tot}(C_{1,1}; -q, Q) =
M(Q_1,q)M(Q_1^{-1},q),\]
to be contrasted with the partition function $Z_{\text{tot}}$ with degree
$d=-2$ obtained earlier in Example~\ref{ex:predef}.
\end{example}
\begin{example}
The following  are the four triangulations corresponding to the crepant
resolutions of the generalised conifold $C_{1,3}:=\{(x,y,z,w)|xy-zw^3=0\}$.
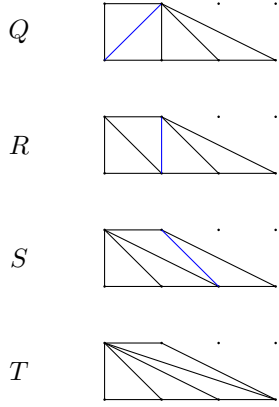
\begin{center}
\begin{tikzpicture}[scale=.75]
\foreach \x in {2,3,4,5} {
  \foreach \y in {-3,-2,-1,0,1,2,3,4} {
    \node at (\x,\y) [circle,fill=black,scale=.1] {};
  }
}
\node at (.5, 3.5) { $Q$ };
\draw (5,3) -- (3,4) -- (2,4) -- (2,3) -- (5,3);
\draw[blue] (2,3) -- (3,4);
\draw (3,3) -- (3,4) -- (4,3);

\node at (.5, 1.5) { $R$ };
\draw (5,1) -- (3,2) -- (2,2) -- (2,1) -- (5,1);
\draw (2,2) -- (3,1);
\draw[blue] (3,1) -- (3,2);
\draw (3,2)--(4,1);

\node at (.5, -.5) { $S$ };
\draw (5,-1) -- (3,0) -- (2,0) -- (2,-1) -- (5,-1);
\draw (3,-1) -- (2,0);
\draw (2,0) -- (4,-1);
\draw[blue] (3,0)--(4,-1);

\node at (.5, -2.5) { $T$ };
\draw (5,-3) -- (3,-2) -- (2,-2) -- (2,-3) -- (5,-3);
\draw (3,-3) -- (2,-2);
\draw (2,-2) -- (4,-3);
\draw (2,-2) -- (5,-3);
\end{tikzpicture}
\end{center}
These  four resolutions  of $xy-zw^3=0$
are obtained from each other by a series of flops,
$$Q \stackrel{\text{\tiny flop} \, \textcolor{blue}{Q_1}} \longrightarrow R
\stackrel{ \text{\tiny flop} \, \textcolor{blue}{R_2}} 
\longrightarrow S \stackrel{\text{\tiny flop} \, \textcolor{blue}{S_3}} \longrightarrow T.$$ 
Let us recall from example \ref{basicflop}
that under a flop of a $(-1,-1)$-line  the formal variable 
changes  from $Q$ to $Q^{-1}$.
Thus the formal variables of the different triangulations are identified as
\begin{equation*} 
\begin{array}{lll}
R_1=Q_1^{-1}, &S_1 = R_1, &  T_1 = S_1, \\
R_2=Q_2,  &S_2 = R_2^{-1}, & T_2 = S_2, \\
R_3=Q_3,  &S_3 = R_3, & T_3 = S_3^{-1},
\end{array}
\end{equation*} 
specifying the isomorphism of second homologies.

Partition functions for the four triangulations are then written down using 
\eqref{prod} as 
\begin{gather*}
\scriptstyle
Z_{top}(q,Q) = M(1,q)^2 \prod_{k=1}^\infty 
(1-Q_1q^k)^{-k} (1-Q_1Q_2q^k)^{-k} (1-Q_1Q_2Q_3q^k)^{-k}
(1-Q_2q^k)^{+k} (1-Q_2Q_3q^k)^{+k} (1-Q_3q^k)^{+k}
\\
\scriptstyle
Z_{top}(q,R) = M(1,q)^2 \prod_{k=1}^\infty (1-R_1q^k)^{-k} (1-R_1R_2q^k)^{+k} (1-R_1R_2R_3q^k)^{+k}
(1-R_2q^k)^{-k} (1-R_2R_3q^k)^{-k} (1-R_3q^k)^{+k}
\\
\scriptstyle
Z_{top}(q,S) = M(1,q)^2 \prod_{k=1}^\infty (1-S_1q^k)^{+k} (1-S_1S_2q^k)^{-k} (1-S_1S_2S_3q^k)^{+k}
(1-S_2q^k)^{-k} (1-S_2S_3q^k)^{+k} (1-S_3q^k)^{-k}
\\
\scriptstyle
Z_{top}(q,T) = M(1,q)^2 \prod_{k=1}^\infty (1-T_1q^k)^{+k} (1-T_1T_2q^k)^{+k} (1-T_1T_2T_3q^k)^{-k}
(1-T_2q^k)^{+k} (1-T_2T_3q^k)^{-k} (1-T_3q^k)^{-k}.
\end{gather*}
Expressing them all in terms of the $Q$ variables, we get
\begin{gather*}
\scriptstyle
Z_{top}(q,Q) = M(1,q)^2 \prod_{k=1}^\infty (1-Q_1q^k)^{-k} (1-Q_1Q_2q^k)^{-k} (1-Q_1Q_2Q_3q^k)^{-k}
(1-Q_2q^k)^{+k} (1-Q_2Q_3q^k)^{+k} (1-Q_3q^k)^{+k}
\\
\scriptstyle
Z_{top}(q,R) = M(1,q)^2 \prod_{k=1}^\infty (1-Q_1^{-1}q^k)^{-k} (1-Q_1^{-1}Q_2q^k)^{+k} (1-Q_1^{-1}Q_2Q_3q^k)^{+k}
(1-Q_2q^k)^{-k} (1-Q_2Q_3q^k)^{-k} (1-Q_3q^k)^{+k}
\\
\scriptstyle
Z_{top}(q,S) = M(1,q)^2 \prod_{k=1}^\infty (1-Q_1^{-1}q^k)^{+k} (1-Q_1^{-1}Q_2^{-1}q^k)^{-k} (1-Q^{-1}_1Q_2^{-1}Q_3q^k)^{+k}
(1-Q_2^{-1}q^k)^{-k} (1-Q_2^{-1}Q_3q^k)^{+k} (1-Q_3q^k)^{-k}
\\
\scriptstyle
Z_{top}(q,T) = M(1,q)^2 \prod_{k=1}^\infty (1-Q_1^{-1}q^k)^{+k} (1-Q_1^{-1}Q_2^{-1}q^k)^{+k} (1-Q_1^{-1}Q_2^{-1}Q_3^{-1}q^k)^{-k}
(1-Q_2^{-1}q^k)^{+k} (1-Q_2^{-1}Q_3^{-1}q^k)^{-k} (1-Q_3^{-1}q^k)^{-k}.
\end{gather*}
Taking the product to assemble the full partition function of the singularity
$C_{1,3}$, after cancellations, we are left with
\begin{equation*} 
\begin{split}
Z^a_{tot}(C_{1,3};q,Q) =
M(1,q)^8 \prod_{k=1}^\infty \scriptstyle
\left(\frac{  (1-Q_1^{-1}q^k) (1-Q_3q^k) (1-Q_1^{-1}Q_2q^k) (1-Q_2^{-1}Q_3q^k)
(1-Q_1^{-1}Q_2Q_3q^k) (1-Q^{-1}_1Q_2^{-1}Q_3q^k)}
{(1-Q_1q^k) (1-Q_3^{-1}q^k) 
(1-Q_1Q_2q^k) (1-Q_2^{-1}Q_3^{-1}q^k) 
(1-Q_1Q_2Q_3q^k) (1-Q_1^{-1}Q_2^{-1}Q_3^{-1}q^k)
}\right)^k,
\end{split}
\end{equation*} 
which can be rewritten, using the expression \eqref{mac} for the generalised 
MacMahon function, as
\begin{equation*} 
\begin{split}
Z^a_{tot} &(C_{1,3};q,Q) =\\
&M(1,q)^8 \frac{M(Q_1^{-1},q) M(Q_1^{-1}Q_2,q) M(Q_1^{-1}Q_2Q_3,q)}
{M(Q_1,q)M(Q_1Q_2,q)M(Q_1Q_2Q_3,q)}
\frac{ M(Q_3,q)M(Q_2^{-1}Q_3,q)  M(Q_1^{-1}Q_2^{-1}Q_3,q) }
{ M(Q_3^{-1},) M(Q_2^{-1}Q_3^{-1},q)M(Q_1^{-1}Q_2^{-1}Q_3^{-1},q)}.
\end{split}
\end{equation*} 
This  expression corresponds to the partition function of a quiver variety
that enjoys a derived equivalence with the crepant resolutions \cite{Na,Y}.
The partition function $Z_{\text{tot}}$ for this case is of vanishing degree,
by \eqref{eq.a}.
\end{example}
To summarise, we defined a partition function for a 
generalised conifold through the 
product of partition functions of all its crepant resolutions. 
The second homologies of the resolutions are identified through a 
canonical ordering of elements, facilitated 
by the absence of homology four-cycles in the resolutions. 
We proved that the new partition function is homogeneous 
with respect to MacMahon factors. This has been contrasted 
with the same product of partition functions with the relation between 
the elements of the second homology group given by Seiberg duality. 

\end{document}